\documentclass[12pt]{article}

\newcommand{\clF}{{\mathcal F}}
\newcommand{\clG}{{\mathcal G}}

\newcommand{\clP}{{\mathcal P}}

\newcommand{\lvG}{{\mathbb G}}

\newcommand{\lvQ}{{\mathbb Q}}

\newcommand{\lvT}{{\mathbb T}}
\newcommand{\lvZ}{{\mathbb Z}}
\newcommand{\Gal}{{\rm Gal}}
\newcommand{\Hom}{\rm Hom}
\newcommand{\Hone}{\textrm{H}^1}
\newcommand{\Inn}{{\rm Inn}}
\newcommand{\hhb}[1]{\hbox to#1pt{}}
\newcommand{\hra}{\hookrightarrow }
\newcommand{\oli}{\overline }
\newcommand{\plim}[1]{\hbox to14pt{lim\kern-14pt\lower4.5pt\hbox{%
$\scriptstyle\longleftarrow$}\kern-8pt\lower8.5pt\hbox{$\scriptstyle{#1}$}}\hhb{3}}
\newcommand{\td}{{\rm td}}
\newcommand{\tlim}{{\tilde\imath} }
\newcommand{\tms}[1]{{#1}^\times\!}
\newcommand{\twoto}{\rightarrow }
\newcommand{\wh}[1]{\widehat{#1}}
\newcommand{\wt}{\widetilde }

\newcommand{\eu}{\mathfrak}

\usepackage{latexsym}
\usepackage{graphics}
\usepackage{amssymb}
\usepackage{amsmath}
\usepackage{amscd}
\usepackage{amsthm}

\newtheorem{thm}{Theorem}

\usepackage{fancyhdr}
\begin{document}

\title{The pro-$p$ Hom-form of the birational anabelian conjecture}

\author{\small {\it Scott Corry} at Appleton and {\it Florian Pop} at Philadelphia}

\date{}

\maketitle

\begin{abstract}
\noindent
We prove a pro-$p$ Hom-form of the birational anabelian
conjecture for function fields over sub-$p$-adic fields. Our
starting point is the corresponding Theorem of Mochizuki
in the case of transcendence degree~1.\footnote{{\it 2000
Mathematics Subject Classification.} Primary 11S20, 12F10,
14G32, 14H30; Secondary 11G20, 14H05.}
\end{abstract}

\section{Introduction}
\label{intro}
\subsection{Grothendieck's anabelian geometry}

On June 27, 1983 A. Grothendieck wrote a letter \cite{G} to
G. Faltings in which he described his dream of an ``anabelian''
algebraic geometry, which has become known as the
Grothendieck Conjecture:
\begin{quote}
{\small A general fundamental idea is that for certain,
so-called ``anabelian'', schemes $X$ (of finite type) over
$K$, the geometry of $X$ is completely determined by
the (profinite) fundamental group $\pi_1(X,\xi) \dots$
together with the extra structure given by the homomorphism:
$$
\pi_1(X,\xi)\rightarrow\pi_1(K,\xi)=\Gal(\oli{K}/K)\qquad (\rm{p}. 280)
$$}
\end{quote}
The intuition that the arithmetic fundamental group of
certain ``anabelian'' schemes over finitely generated
fields, $K$, should be so extraordinarily rich as to completely
specify the isomorphism type of the scheme has been
borne out over the last quarter of a century, especially in
the case of hyperbolic curves, which were the main class
for which Grothendieck claimed the anabelian title 
\cite{M}, \cite{N1}, \cite{N2}, \cite{T}. But Grothendieck
also expected the function fields of algebraic $K$-varieties
to be anabelian, giving rise to the birational version of the
Grothendieck Conjecture which can be viewed as a
generalization of the Galois characterization of global
fields due to J. Neukirch, K. Uchida, et. al. \cite{Neu},
\cite{U}, \cite{P1}, \cite{P2}.

In the years since Grothendieck's letter, the anabelian
philosophy has evolved to include various expectations
about recovering schemes from quotients of the arithmetic
fundamental group, such as the pro-$p$ situation
considered here. Furthermore, it has become clear that
one need not restrict oneself to finitely generated base
fields, and in fact anabelian phenomena even occur over
algebraically closed base fields, i.e. in the absence of any
arithmetic structure \cite{P4}. Perhaps most importantly,
Sh.~Mochizuki discovered in \cite{M} that
\begin{quote}
{\small \dots the Grothendieck Conjecture for hyperbolic
curves [is] an essentially local, $p$-adic result that belongs
to that branch of arithmetic geometry known as $p$-adic
Hodge theory. (p. 326)} 
\end{quote}

Following Mochizuki, we are concerned here with the
question of recovering embeddings of function fields over
sub-$p$-adic fields from the induced homomorphisms of
pro-$p$ Galois groups.

\subsection{Statement of the Theorem}
Let $k$ be a fixed sub-$p$-adic field, i.e., a subfield of a
function field $k_1|\lvQ_p$. Fix $\oli k$, an algebraic
closure of $k$, and denote by $G_k$ the absolute Galois
group. Let $\clF_k$ be the category of regular function
fields $K|k$, and $k$-embeddings of function fields.
Further, let $\clG_k$ be the category of profinite groups,
$G$, endowed with a surjective augmentation morphism,
$\pi_G:G\to G_k$, such that $\ker(\pi_G)$ is a pro-$p$
group, and outer open $G_k$-homomorphisms, i.e.,
a morphism from $G$ to $H$ in $\clG_k$ is of the form
$\Inn_{G_k}(H)\circ f$, where $f:G\to H$ is an open
homomorphism such that $\pi_G=\pi_H\circ f$, and
$\Inn_{G_k}(H)$ denotes the group of inner
automorphisms of $H$ which lie over $G_k$.
Since $G_k$ has trivial center, $\Inn_{G_k}(H)$ consists
exactly of the inner automorphisms of $H$ defined by an
element of $\ker(\pi_H)$. Finally, we remark that there
exists a naturally defined functor from $\clF_k$ to $\clG_k$: 
for $K|k$ from $\clF_k$, let $\wt K|K\oli k$
be a maximal pro-$p$ extension. Then $\wt K|K$ is
Galois, and $\Pi_K:= \Gal(\wt K|K)$ endowed
with the projection $\pi_K:\Pi_K\to G_k$ is an object
of $\clG_k$. Further, a morphism $\imath:K|k\hra L|k$
in $\clF_k$ extends uniquely to a $\oli k$-embedding
$\oli\imath:K\oli k\hra L\oli k$, and $\oli\imath$ has
prolongations $\tlim:\wt K\hra\wt L$. Each
such prolongation $\tlim:\wt K\hra\wt L$ gives
rise to an open $G_k$-homomorphism
$\Phi_{\tlim}:\Pi_L\to\Pi_K$ defined by
$$
\Phi_{\tlim}(g)=\tlim^{-1}\circ g \circ \tlim\,,
  \quad g\in\Pi_L\,,
$$
and any two such prolongations are conjugate by
an element from $\ker(\pi_K)$. Thus, sending
each $K|k$ from $\clF_k$ to $\pi_K:\Pi_K\to G_k$ in $\clG_k$
yields a well defined functor from $\clF_k$ to $\clG_k$.


The purpose of this note is to prove the following
{\it Galois by pro-$p$ Hom-form\/} of the birational
anabelian conjecture:

\begin{thm}\label{main}
The above functor from $\clF_k$ to $\clG_k$ is fully
faithful, i.e., for regular function fields $K|k$ and $L|k$,
there is a canonical bijection
$$
{\rm Hom}_k(K,L) \to \Hom_{\clG_k}(\Pi_L,\Pi_K).
$$
Equivalently, for fixed field extensions $\wt K|K$ and
$\wt L|L$ as above, the map
$$
\tlim\mapsto\Phi_{\tlim}\ \ \ \hbox{\it with \ }
  \Phi_{\tlim}(g)=\tlim^{-1}\circ g \circ \tlim
    \hbox{\it \ \ for \ }  g\in\Pi_L\,,
$$
is a bijection from the set of $\oli k$-embeddings
$\tlim:\wt K\hra\wt L$ onto the set of all the open
$G_k$-morphisms $\Pi_L\to\Pi_K$.
\end{thm}

Before embarking on the proof, the following
comments are in order: first, if $\td(K|k)=1$, then
Theorem~\ref{main} is a special case of Theorem 16.5
in the fundamental paper by Mochizuki~\cite{M}. Second,
the above Theorem~\ref{main} implies the corresponding
{\it full profinite version,\/} in which $\Pi_K$ is replaced
by the full absolute Galois group $G_K$ of $K$. But naturally,
the above Theorem~\ref{main} does not follow from the
corresponding full profinite version. Finally, the full
profinite version of Theorem~\ref{main} above was
proved by Mochizuki in loc.~cit., where it appears as
Corollary 17.1. There he uses an inductive procedure on
$\td(K|k)$ which is ill-suited to the pro-$p$ situation, and
hence he obtains only a profinite result. In our proof of
Theorem~\ref{main} we use Mochizuki's pro-$p$ result
for the transcendence degree one case, but instead of
proceeding inductively on the transcendence degree,
we will make use of the second author's ideas as described
in~\cite{P3}.


\section{Proof of Theorem \ref{main}}
\label{sec:1}
The proof of Theorem~\ref{main} will have two parts:
\begin{enumerate}
\item[i)] Given an open $G_k$-homomorphism
$\Phi:\Pi_L\to\Pi_K$, there exists a $\oli k$-embedding
$\tlim_\Phi:\wt K\hra\wt L$ which defines $\Phi$ as
indicated above, i.e., such that $\Phi=\Phi_{\tlim_\Phi}$.

\vskip3pt

\item[ii)] The map $\tlim\mapsto\Phi_{\tlim}$
is injective.
\end{enumerate}

First, let us recall the following basic facts
about $p$-adic completions of abelian groups (which
should not be confused with the pro-$p$ completions
of such groups). For every abelian group $A$, let
$\jmath_A:A\to\wh A:=\plim{e}\, A/p^e$ be\break
\vskip-13pt\noindent
the $p$-adic completion homomorphism from $A$
to its $p$-adic completion $\wh A$. It is clear that the
passage from $A$ to $\wh A$ is a functor from the category
of abelian groups (which is the same as the category
of $\lvZ$-modules) to the category of $\lvZ_p$-modules.
Moreover, ${\rm ker}(\jmath_A)=A_\infty:=\{\,a\in A \mid
   \forall e>0\ \exists b\in A \ \hbox{s.t.} \ a=p^eb\,\}$
is the maximal $p^\infty$-divisible subgroup of $A$. In
particular, if $0\to A'\to A$ is a short exact sequence
of abelian groups, then the resulting canonical sequence
$0\to\wh A'\to \wh A$ is exact if and only if
$A'_\infty=A'\cap A_\infty$. Hence if $A/A'$ has trivial
$p$-torsion, then $0\to\wh A'\to \wh A$ is exact. In the
sequel, we will consider/use the $p$-adic completions of
the multiplicative groups $k^\times\subset K^\times$
of field extensions $K|k$ from $\clF_k$, which we denote
simply by $\wh k:=\wh{k^\times}$, $\wh K:=\wh{K^\times}$,
etc.

We remark that for
$K|k$ in $\clF_k$ one has:
\vskip3pt

a) $K^\times\!/k^\times$ is a free abelian group,
hence $\wh k\hra \wh K$.
\vskip3pt

Indeed, this follows from the observation that $K^\times/k^\times$ is
the subgroup of principal divisors inside the free abelian group ${\rm Div}(X)$,
where $X$ is any projective normal model of the function
field $K|k$.
\vskip3pt

b) ${\rm ker}(\jmath_K)=\mu'$ is the finite group of roots
of $1$ of order prime to $p$ in $k$.
\vskip3pt

Indeed, recall that ${\rm ker}(\jmath_K)$ equals the
$p^\infty$-divisible subgroup of $K^\times$. Now since
$K^\times\!/k^\times$ is a free abelian group, it follows
that ${\rm ker}(\jmath_K)\subset k^\times$. And since
$k$ is a sub-$p$-adic field, it is embeddable into a
regular function field $k_1|k_0$, where $k_0|\lvQ_p$
is a finite field extension. But then reasoning as above,
$k_1^\times\!/k_0^\times$ is a free abelian group, and
so ${\rm ker}(\jmath_K)\subset k_0^\times$. Hence
finally, ${\rm ker}(\jmath_K)$ is contained in the
$p^\infty$-divisible subgroup of $k_0^\times$, which is the
group $\mu'_{k_0}$ of roots of unity of order prime
to $p$ in $k_0$, hence a finite group. Thus
${\rm ker}(\jmath_K)=\mu'_{k_0}\cap k=:\mu'$.

\vskip5pt

We next recall the following basic facts from
Kummer Theory:  For every $K|k$ from $\clF_k$, let
$\lvT_{p,k}=\plim{e}\,\mu_{p^e\!,\, \oli k}$ and
$\lvT_{p,K}=\plim{e}\,\mu_{p^e\!,\, \wt K}$
be the Tate modules of $\lvG_{m,k}$, respectively
$\lvG_{m,K}$. Then via the canonical inclusion
$\oli\imath:\oli k\to\wt K$ we can/will identify $\lvT_{p,K}$
with $\lvT_p:=\lvT_{p,k}$. Then Kummer Theory yields
a canonical isomorphism of $p$-adically complete groups
\[
\delta_K:\wh K=\plim{e}\, K^\times\!/p^e\to
 \plim{e}\,\Hone(\Pi_K,\mu_{p^e\!\!,\wt K})=\Hone(\Pi_K,\lvT_p).
\]
Therefore we will make the identification
$\wh K=\Hone(\Pi_K,\lvT_p)$, if this does not lead
to confusion. By the functoriality of Kummer
Theory, the surjective projection $\pi_K:\Pi_K\to G_k$
gives rise to a canonical homomorphism
$
\Hone(\pi_K):\wh k\hra \wh K
$
which is nothing but the $p$-adic completion of the
structural morphism $k\hra K$, and it is an embedding
by remark~a) above.
Furthermore, if $K|k$ and
$L|k$ are objects from $\clF_k$, and $\Phi:\Pi_L\to\Pi_K$
is an open $G_k$-morphism, then by functoriality
we get an embedding of $p$-adically
complete groups
$$
\Hone(\Phi):\wh{K}=\Hone(\Pi_K,\lvT_p)\to
    \Hone(\Pi_L,\lvT_p)=\wh{L}
$$
which identifies $\wh k\subset\wh K$ with
$\wh k\subset\wh L$. Finally note that if
$\Phi=\Phi_{\tlim}$ is defined by a morphism
$\imath:K|k\hra L|k$ from $\clF_k$, then
$\Hone(\Phi_{\tlim})=\wh\imath$ is nothing but
the $p$-adic completion of the $k$-embedding
$\imath:K|k\hra L|k$, and therefore one has:
$$
\Hone(\Phi_{\tlim})\circ\jmath_K(x)=
      \jmath_L\circ\imath(x),\quad x\in K^\times.
  \leqno{\indent(\dag)}
$$


\noindent
\textbf{Proof of i): }

\vskip3pt
\noindent
{\it Claim 1.\/}  $\Hone(\Phi)\circ\jmath_K(K^\times)
   \subseteq\jmath_L(L^\times)$.

\vskip3pt
\noindent
{\it Proof of Claim 1\/}: Consider $t\in K^\times$
arbitrary. First, if $t\in k^\times$, then $\Hone(\Phi)$
identifies $\jmath_K(t)$ with $\jmath_L(t)$ by the
discussion above. Second, let
$t\in K^\times\backslash k^\times$. Then the
inclusion $k(t)\subseteq K$ is a morphism in
$\clF_k$, hence gives rise canonically to a
$\clG_k$-morphism $\Phi_t:\Pi_K\to\Pi_{k(t)}$. But then
$\Phi_t\circ\Phi:\Pi_L\to\Pi_{k(t)}$ is a
$\clG_k$-morphism too. Since $\td\big(k(t)|k\big)=1$,
it follows by Theorem~16.5 from \cite{M} that
$\Phi_t\circ\Phi$ is defined by a
$k$-embedding $\imath_t:k(t)\to L$, i.e., is of the
form $\Phi_t\circ\Phi=\Phi_{\tlim_t}$. Hence
by the assertion~$(\dag)$ above,
$\Hone(\Phi_t\circ\Phi):\wh{k(t)}\to\wh{L}$ is
exactly the $p$-adic completion of $\imath_t:k(t)\to L$,
and we get:
$$
\Hone(\Phi_t\circ\Phi)\circ\jmath_{k(t)}(x)=
            \jmath_L\circ\imath_t(x)\in \jmath_L(L^\times),
                               \quad x\in k(t)^\times.
$$
By functoriality,
$\Hone(\Phi_t\circ\Phi)=\Hone(\Phi)\circ\Hone(\Phi_t)$,
and $\Hone(\Phi_t)$ is the $p$-adic completion of the
inclusion $k(t)\hra K$.
Hence $\Hone(\Phi_t)\circ\jmath_{k(t)}(x)=\jmath_K(x)$
for $x\in k(t)^\times$. Combining these equalities, we
finally get
$$
\Hone(\Phi)\circ\jmath_K(x)=
      \jmath_L\circ\imath_t(x)
\in \jmath_L(L^\times),
                               \quad x\in k(t)^\times.
\leqno{\indent(\dag')}
$$
Since $t$ was arbitrary, 
this concludes the proof of Claim~1. $\Box$


Next let us identify $\tms K/\mu'$ and $\tms L/\mu'$
with their images in $\wh{K}$, respectively $\wh{L}$,
via the $p$-adic completion homomorphisms $\jmath_K$,
respectively $\jmath_L$. Then by Claim~1 above,
$\Hone(\Phi)$ maps $\tms K/\mu'$ into $\tms L/\mu'$,
and identifies $\tms k/\mu'\subset \tms K/\mu'$ with
$\tms k/\mu'\subset \tms L/\mu'$. Further, $\Hone(\Phi)$
is nothing but the $p$-adic completion of its restriction
to $\tms K/\mu'$. Modding out by $\tms k/\mu'$ thus
yields an embedding of free abelian groups
$\alpha:\tms K/k^\times\hra\tms L/k^\times$
canonically defined by $\Hone(\Phi)$.
Now we regard $(K,+)$ and $(L,+)$ as infinite dimensional
$k$-vector spaces, and denote by  $\clP(K):=\tms K/k^\times$ and
$\clP(L):=\tms L/k^\times$ their projectivizations.
Then $\alpha:\clP(K)\hra\clP(L)$  is an inclusion which
respects the multiplicative structures.

\vskip5pt
\noindent
{\it Claim 2.\/} The map $\alpha:\clP(K)\hra\clP(L)$ preserves
lines.

\vskip3pt
\noindent
{\it Proof of Claim 2\/}:
A line in $\clP(K)$ is the image of a two-dimensional
$k$-subspace of $K$, say ${\eu l}_{t_1t_2}:=k t_1+k t_2$,
where $t_1, t_2 \in K$ are $k$-linearly independent.
Note that ${\eu l}_{t_1t_2}=t_1\cdot{\eu l}_t$,
where $t=t_2/t_1$ and ${\eu l}_t:=k+k t$. Since
multiplication by $t_1/k^\times$ is a line-preserving automorphism
of $\tms K/k^\times=\clP(K)$,
$\alpha$ is multiplicative, and multiplication by
$\alpha(t_1/k^\times)$ is a line-preserving automorphism of $\clP(L)$, it 
suffices to show
that $\alpha$ maps the lines ${\eu l}_t$, $t\in K\backslash k$, 
to lines in $\clP(L)$. In order to do this, we need only remark
that by relation~$(\dag')$ above we have:
$$
\Hone(\Phi)\circ\jmath_K(kt+k)^\times=
      \jmath_L\circ\imath_t(kt+k)^\times=
           \jmath_L\big(k\imath_t(t)+k\big)^\times.
$$
Thus ${\eu l}_t\subset\clP(K)$ is mapped bijectively
onto ${\eu l}_{\imath_t(t)}\subset\clP(L)$. $\Box$

\vskip3pt

Let $L^\times_0\subseteq L^\times$ denote the
preimage of $\alpha\big(\clP(K)\big)$ in $L$, and
let us set $L_0:=L_0^\times\cup\{0_L\}\subseteq L$.
Since $\alpha$ preserves lines, it follows that $L_0$
is a field containing $k$, and $\alpha:\clP(K)\to\clP(L_0)$
is a line-preserving bijection. By the Fundamental
Theorem of projective geometry (Theorem 2.26 of
\cite{A}),  we conclude that $\alpha$ is induced by a
$k$-semilinear isomorphism of $k$-vector spaces
$\alpha_K:(K,+)\to(L_0,+)$, which is unique up to
$k$-semilinear homotheties. 

\vskip3pt
\noindent
{\it Claim 3.\/} $\alpha_K$ is $k$-linear.

\vskip3pt
\noindent
{\it Proof of Claim 3\/}:
Let $\mu:k\to k$ be the field isomorphism with
respect to which $\alpha_K$ is $k$-semilinear, i.e.
$\alpha_K(ax)=\mu(a)\alpha_K(x)$ for all $a\in k, x\in K$.
We wish to show that $\mu=\textrm{id}_k$. For this,
pick $t\in K^\times\backslash k^\times$ and consider
the inclusion of $k$-vector spaces $(k(t),+)\subseteq (K,+)$.
Then $\alpha_t:=\alpha_K|_{k(t)}$ is a $k$-semilinear
map with respect to $\mu$, with  projectivization
$$
\clP(\alpha_t)=\clP(\alpha_K)|_{\clP(k(t))}=
                               \alpha|_{\clP(k(t))}=\clP(i_t),
$$
the last equality coming from relation~$(\dag')$
interpreted in terms of projectivizations. This immediately
implies that $\text{im}(i_t)\subseteq L_0$, so $i_t:k(t)\to L_0$
is a $k$-embedding a fields, \emph{a fortiori} a $k$-linear
map of $k$-vector spaces. By uniqueness, $i_t$ differs from
$\alpha_t$ by a $k$-semilinear homothety, and in particular,
$i_t$ must be $k$-semilinear with respect to $\mu$, so
that that $\mu=\textrm{id}_k$ as claimed.~$\Box$

\vskip3pt

As in \cite{P3}, it now follows that
setting $\imath_K:=\alpha_K(1_K)^{-1}\cdot\alpha_K$,
the resulting map $\imath_K:K\to L$ is actually a
$k$-isomorphism of fields, whose projectivization equals
$\alpha$. In particular, the $p$-adic completion of $\imath_K$
equals $\Hone(\Phi)$, and $\imath_K$ is the unique
embedding of fields $K\hra L$ with this property.


Now let $K'|K$ be a finite Galois extension contained
in $\wt K$, and $k':=K'\cap\oli k$. Then $\Pi_{K'}$
is an open normal subgroup of $\Pi_K$, and
$\Phi^{-1}(\Pi_{K'})=\Pi_{L'}$ for some finite Galois
extension $L'|L$ contained in $\wt L$. Since $\Phi$ is
a $G_k$-homomorphism, it follows that $L'\cap\oli k=k'$,
and $K'|k'$ and $L'|k'$ are regular function fields over
$k'$ which is a sub-$p$-adic field. The restriction $\Phi'$
of $\Phi$ to $\Pi_{L'}$ yields an open 
$G_{k'}$-homomorphism $\Phi':\Pi_{L'}\twoto\Pi_{K'}$.
{\it Mutatis mutandis,\/} we obtain from $\Phi'$ a
$k'$-embedding $\imath_{K'}:K'\hra L'$ such that its
$p$-adic completion is $\Hone(\Phi')$. The compatibility
relation
$
{\rm res}_{\Pi_L}^{\Pi_{L'}}\circ\Hone(\Phi')=
    \Hone(\Phi)\circ{\rm res}_{\Pi_K}^{\Pi_{K'}}
$
translates into the fact that $\imath_K$ is the restriction
of $\imath_{K'}$ to $K$. Note that extensions of the 
form $K'|K$ above exhaust $\wt K|K$, so
taking limits we obtain a $\oli k$-embedding
$\tlim_\Phi:\wt K\hra\wt L$ such that
$(\tlim_\Phi)|_{K'}=\imath_{K'}$ for all $K'|K$.

\vskip3pt
\noindent
{\it Claim 4.\/}  $\Phi_{\tlim_\Phi}=\Phi$.

\vskip3pt
\noindent
{\it Proof of Claim 4\/}:
Indeed, for $K'|K$ and the corresponding $L'|L$ as above,
$\Phi$ yields a surjection
$
\Psi':\Gal(L'|L)\twoto\Gal(K'|K)\,.
$
Note that $\Hone(\Pi_{K'},\lvT_p)$ is a $\Gal(K'|K)$
module, and correspondingly for $L'|L$. Moreover, for all
$\sigma\in\Gal(L'|L)$ we have the following commutative
diagram:
$$
\begin{CD}
\Hone(\Pi_{K'},\lvT_p)@>\Hone(\Phi')>>\Hone(\Pi_{L'},\lvT_p)\\
@V\Psi'(\sigma)VV                          @VV\sigma V\\
\Hone(\Pi_{K'},\lvT_p)@>\Hone(\Phi')>>\Hone(\Pi_{L'},\lvT_p).
\end{CD}
$$
which via the Kummer Theory isomorphisms
translates into:
$$
\imath_{K'}\circ\Psi'(\sigma)=\sigma\circ\imath_{K'},
\quad \sigma\in\Gal(K'|K).
$$
Hence taking limits over all the $K'|K$, we get as required:
$$\tlim_\Phi\circ\Phi(g)=g\circ\tlim_\Phi,
\quad g\in\Pi_L\,. \ \Box
$$

\noindent
\textbf{Proof of ii): } Let $\tlim:\wt K\hra \wt L$ be a
$\oli k$-embedding, and $\Phi_{\tlim}:\Pi_L\to\Pi_K$
be the corresponding
open $G_k$-homomorphism.  For every finite Galois
sub-extension $L'|L$ of $\wt L|L$ and the corresponding
$K':=\tlim^{-1}(L')$, the restriction of $\Phi_{\tlim}$
to $\Pi_{L'}$ defines an open $G_{k'}$-homomorphism
$\Phi_{\tlim,L'}:\Pi_{L'}\to\Pi_{K'}$, where $k':=\oli k\cap L'$.
By the discussions above, the Kummer morphism
$\Hone(\Phi_{\tlim,L'})$ is the $p$-adic completion
of $\imath_{K'}:=\tlim|_{K'}$, and in particular,
$\Hone(\Phi_{\tlim,L'})$ determines $\imath_{K'}$
uniquely. By taking limits, it follows that $\tlim$
is uniquely determined by the family of Kummer
morphisms $\Hone(\Phi_{\tlim,L'})$ with $L'|L$
as above. Assertion ii) is thus proven. $\Box$

\vskip3pt
\noindent
This concludes the proof of Theorem~\ref{main}.

\center{---------------------------------}

\center{\small Lawrence University, Dept. of Mathematics,
P.O. Box 599, Appleton, WI 54912, USA \\
e-mail: scott.corry@lawrence.edu}
\center{\small University of Pennsylvania, Dept. of Mathematics,
209 S. 33rd St., Philadelphia, PA 19104, USA \\
e-mail: pop@math.upenn.edu}

\end{document}